\documentstyle[11pt,titlepage]{report}
\begin{document}
\title {
\vspace*{-7ex}
{
\begin{flushright}
	  \textbf{\large LANL xxx E-print archive No. dg-ga/9704004}\\[5ex]
\end{flushright}
}
	    \hspace{0.7cm}     \bf\huge {Transports along Paths} \\
   	    \hspace{1.7cm}	\bf\huge {in Fibre Bundles} 	 \\
							\vspace{0.9cm}
	  {\bf\LARGE III. Consistency with Bundle Morphisms}  }
\vspace {28pt}
\renewcommand{\thefootnote}{\fnsymbol{footnote}}
\author{Bozhidar Zakhariev Iliev
\thanks {\sl Permanent address:
Institute for nuclear research and nuclear energy, Bulgarian
academy of sciences, blvd. Tzarigradsko
chaus\`ee 72, 1784 Sofia, Bulgaria.}
\thanks{E-mail address: bozho@inrne.bas.bg}
\thanks{URL: http://www.inrne.bas.bg/mathmod/bozhome/}
}
\date {\vspace{4ex}
Published: Communication JINR, E5-94-41, Dubna,  1994\\
LANL xxx archive E-print No. dg-ga/9704004}
\maketitle

\def\thesection {\arabic{section}.}

\setcounter {chapter} {1}
\setcounter {equation} {0}
\section {\bf INTRODUCTION}
\par
In Ref. [1] the problem has been investigated on consistency, or
compatibility, of linear transports along paths in vector bundles
and bundle metrics, i.e. when the transports preserve the (scalar
products defined by the) metric. The present paper generalizes this
problem to and deals with the problem of consistency (compatibility) of
arbitrary transports along paths in fibre bundles [2] and acting
between these fibre bundles bundle morphisms [3,4]. This task is
sufficiently general, to cover from a unified point of view, all
analogous problems from the  literature available to the author.
\par
The problem for consistency of transports along paths and
bundle morphisms is stated in a general form in Sect.~2. Some necessary
and sufficient conditions for such consistency are found. It
is proved that the introduced concept for consistency is a special
case of the one for sections of a fibre bundle transported by means
of a transport along paths [2]. In Sect.~3 these concept and results
are applied to the special case of a transport along paths
and bundle morphisms acting in one and the same fibre bundle. Also
some examples are presented. Sect.~4 contains a detailed investigation
of the consistency between leaner transports along paths in a
vector bundle and a Hermitian structure in it. Sect.~5 closes the
paper with some concluding remarks.
\par
Below we summarize certain definitions
and results from [2] needed for this paper.
	\par
	By $(E,\pi ,B)$ is denoted an arbitrary (topological) fibre bundle
with a base $B$, bundle space $E$, projection $\pi :E\to B$, and homeomorphic
fibres $\pi ^{-1}(x), \  x\in B$ [4-6].
	\par
The set of sections of $(E,\pi ,B)$ is $Sec(E,\pi ,B)$, i.e.
$\sigma \in Sec(E,\pi ,B)$ means $\sigma :B \to E$ and
$\pi \circ \sigma ={\it id}_B$, where ${\it id}_{X}$ is the identity
map of the set X.
	\par
By $J$ and $\gamma :J \to B$ are denoted, respectively, an arbitrary real
interval and a path in $B$.
	\par
The transport along paths in $(E,\pi ,B)$ is a map
$I:\gamma \mapsto I^{\gamma }$, where
$I^{\gamma }:(s,t)\mapsto I^{\gamma }_{s\to t}, \  s,t\in J$ in which the maps
$I^{\gamma }_{s\to t}:\pi ^{-1}(\gamma (s)) \to \pi ^{-1}(\gamma (t))$
satisfy the equalities
\begin{eqnarray}
& I^{\gamma }_{t\to r}\circ I^{\gamma }_{s\to t}=I^{\gamma }_{s\to r},
\quad r,s,t\in J, & \\
& I^{\gamma }_{s\to s}={\it id}_{\pi ^{-1}(\gamma (s))}, \quad s\in J, &
\end{eqnarray}
and its general form is described by
\begin{equation}
 I^{\gamma }_{s\to t}={\left( F^{\gamma }_{t} \right)}^{-1}\circ
 F^{\gamma }_{s}, \quad s,t\in J,
\end{equation}
$F^{\gamma }_{s}:\pi ^{-1}(\gamma (s)) \to Q, \  s\in J$ being one-to-one maps
onto one and the same set Q.
	\par
 In the case of a linear transport along paths in a vector bundle [7]
the corresponding to (1.3) general form of the transport matrix is (see
[7], proposition 2.4)

\begin{equation}
 H(t,s;\gamma )=F^{-1}(t;\gamma )F(s;\gamma ), \quad s,t\in J
\end{equation}
in which $F(s;\gamma )$ is a nondegenerate matrix function.

\setcounter {chapter} {2}
\setcounter {equation} {0}
\section {\bf GENERAL THEORY}
	Let there be given two fibre bundles
$\xi _{h}:=(E_{h},\pi _{h},B_{h}),\quad h=1,2$ in
which defined are, respectively, the transports along paths ${}^1 \! I$ and
${}^2 \! I$. Let $(F,f)$ be a bundle morphism from $\xi _{1}$ into $\xi _{2}$,
i.e. (see [3,4])
$F:E_{1} \to E_{2}, \quad f:B_{1} \to B_{2}$ and
$\pi _{2}\circ F=f\circ \pi _{1}$.
Let $F_{x}:=F\mid \pi ^{-1}_{1}(x)$ for $x\in B_{1}$ and $\gamma :J \to
B_{1}$ be an arbitrary path in $B_{1}$.
	\par
 	{\bf Definition 2.1.} The bundle
morphism $(F,f)$ and the pair $({}^1 I,{}^2I)$ of transports, or the
transports ${}^1\! I$ and ${}^2\! I$, along paths will be called consistent
(resp. along the path $\gamma $) if they commute in a sense that the
equality
\begin{equation}
F_{\gamma (t)}\circ {}^1 I^{\gamma }_{s\to t}=
{}^{2} I^{f\circ \gamma}_{s \to t}\circ F_{\gamma (s)}, \quad s,t\in J
\end{equation}
is fulfilled for every (resp. the given) path $\gamma $.
	\par
This definition contains as an evident special case definition 1.1 from [1].
In fact, to prove this it is sufficient to put in it: $\xi _{1}=(E,\pi
,B)\times (E,\pi ,B)$, where $(E,\pi ,B)$ is a vector bundle;
$\xi _{2}=({\bf R},\pi _{0},0)$, where $0\in {\bf R}$ and
 $\pi _{0}:{\bf R} \to \{0\}$; $F_{x}=g_{x}$, where $g_{x}:\pi ^{-1}(x)\times
\pi ^{-1}(x) \to {\bf R}, \quad x\in B$ are nondegenerate symmetric and
bilinear maps; $f:B\times B \to \{0\}\subset {\bf R}$;
${}^{1}\! I^{\gamma }_{s\to t}=I^{\gamma }_{s\to t}\times
I^{\gamma }_{s\to t}$, where $I^{\gamma }$ is a transport along
 $\gamma :J \to B$ in $(E,\pi ,B)$ and $s,t\in J$;
${}^{2}\! I^{f\circ \gamma }_{s \to t}={\it id}_{\footnotesize\bf R}$.
 	\par
 Analogously to the considerations in [1], Sect.~2 here can be
formulated in a general form, for example, the following problems: to be
found necessary and/or sufficient conditions for consistency between bundle
morphisms and (ordered) pairs of transports along paths; to be found, if
any, all bundle morphisms (resp. transports along paths) which are
consistent with a given pair of transports along paths (resp. bundle
morphism), etc. Further we will consider some results in this field.
	\par
 	{\bf Proposition 2.1.} The bundle morphism $(F,f)$ and the pair
 $({}^{1}\! I,{}^{2}\! I)$
of transports along paths are consistent (resp. along the path $\gamma )$
iff there exist $s_{0}\in J$ and a map
\begin{equation}
 C(s_{0};\gamma ,f\circ \gamma ):\pi ^{-1}_{1}(\gamma (s_{0})) \to
 \pi ^{-1}_{2}((f\circ \gamma )(s_{0})),
\end{equation}
such that
\begin{equation}
 F_{\gamma (s)}={}^{2}\! I^{f\circ \gamma }_{s_{0}\to s}\circ
 C(s_{0};\gamma ,f\circ \gamma )\circ {}^{1}\! I^{\gamma }_{s\to s_{0}}
\end{equation}
for every (resp. the given) path $\gamma $.
	\par
 	{\bf Remark.} The reason for which as arguments of $C$ are written
$s_{0}$, $\gamma $ and $f\circ \gamma $ will be cleared up below in
proposition 2.2, where its general structure is described.
	\par
	{\bf Proof.}
 If $(F,f)$ is consistent with $({}^{1}\! I,{}^{2}\! I) $ (resp. along $\gamma
$), then, by definition, (resp. along $\gamma $) is valid (2.1), which, due
to $\left( {I^{\gamma }_{s\to t}} \right) ^{-1}=I^{\gamma }_{t\to s}$ (see
(1.3)), is equivalent to

	\begin{equation}
 F_{\gamma (t)}={}^{2}\! I^{f\circ \gamma}_{s \to t}\circ
 F_{\gamma (s)}\circ {}^{1}\! I^{\gamma }_{t\to s}, \quad s,t\in J.
 	\end{equation}
Fixing $s_{0}\in J$
and putting $s=s_{0}$ and $t=s$  in this equality, we get (2.3) (resp. along
$\gamma $) with $C(s;\gamma ,f\circ \gamma )=F_{\gamma (s_{0})}$.
 And vice versa, if (2.3) is true (resp. along $\gamma $) for some
$s_{0}$ and $C$, then due
to (1.1), we have
\begin{eqnarray}
 & & F_{\gamma (t)}\circ{}^{1}\! I^{\gamma }_{s\to t} =
{}^{2}\! I^{f\circ \gamma }_{s_{0} \to t}\circ
 C(s_{0};\gamma ,f\circ \gamma )\circ {}^{1}\! I^{\gamma }_{t\to s_{0}}
\circ {}^{1}\! I^{\gamma }_{s\to t}= \nonumber \\
 & & ={}^{2}\! I^{f\circ \gamma }_{s\to t}\circ {}^{2}\! I^{f\circ
\gamma}_{s_{0}\to s} \circ C(s_{0};\gamma ,f\circ \gamma ) \circ {}^{1}\!
I^{\gamma }_{s\to s_{0}}={}^{2}\! I^{f\gamma }_{s\to t} \circ F_{\gamma (s)},
\nonumber \end{eqnarray} i.e. (2.1) is identically satisfied (resp. along
$\gamma $) for $s,t\in J$, so $(F,f)$ and $({}^{1}\! I,{}^{2}\! I)$ are consistent
(resp. along $\gamma $).\hbox{\vrule width 6pt height 6pt depth 0pt}
 	\par
	{\bf Corollary 2.1.} The equality (2.3) is valid iff
		\setcounter{equation}{2}
		\def\theequation{\thechapter.\arabic{equation}'}
	\begin{equation}
F_{\gamma (s)}={}^{2}\! I^{f\circ \gamma }_{t_{0} \to s}\circ
 C(t_{0};\gamma ,f\circ \gamma )\circ {}^{1}\! I^{\gamma }_{s \to t_{0}},
\end{equation}
	\def\theequation{\thechapter.\arabic{equation}}
	\setcounter{equation}{4}
where
	\begin{equation}
C(t_{0};\gamma ,f\circ \gamma )={}^{2}\! I^{f\circ \gamma }_{s_{0} \to
t_{0}}\circ C(s_{0};\gamma ,f\circ \gamma )\circ {}^{1}\! I^{\gamma
}_{t_{0}\to s_{0}} \end{equation} for arbitrary $t_{0}\in J$, i.e. the
existence of $s_{0}\in J$ and a map (2.2) for which (2.3) is valid leads to
the existence of maps (2.5) for which (2.3$'$) is true for every $t_{0}\in J$
and vice versa.  \par We have written (2.5) as an equality, but not as a
definition of $C(t_{0};\gamma ,f\circ \gamma )$, because from (1.3) and (2.5),
considered as a definition of $C$, for arbitrary $t_{0}\in J$ and fixed
$s_{0}\in J$ the validity of (2.5) follows for every $t_{0},s_{0}\in $J.  \par
{\bf Proof.} If (2.3) is valid, then by proposition 2.1 $(F,f)$ and $({}^{1}\!
	I,{}^{2}\!  I)$ are consistent, i.e. (2.4) is true, from where (2.3$'$)
with $C(t_{0};\gamma ,f\circ \gamma )= F_{\gamma (t_{0})}$ follows.  In the
opposite direction the proposition follows from the substitution of (2.5) into
(2.3$'$) and the usage of (1.1). (The same result follows also directly from
(2.3) and ${\left( I^{\gamma}_{s\to t} \right)}^{-1}=I^{\gamma }_{t\to s}$
(see (1.3))).\hbox{\vrule width 6pt height 6pt depth 0pt}
	\par {\bf Corollary 2.2.} The bundle morphism $(F,f)$ and
the pair $({}^{1}\! I,{}^{2}\! I)$ of transports along paths are consistent (resp. along
the path $\gamma $) iff there exists a bundle morphism $(C,f)$ from $\xi
_{1}$ on $\xi _{2}$ such that for every (resp. the given)
$\gamma :J \to B_{1}$ is fulfilled:
a) $(C,f)\mid \gamma (J):=(C\mid \pi ^{-1}_{1}(\gamma (J)),f\mid \gamma (J))$
 is a bundle morphism from $\xi _{1}\mid _{\gamma (J)}$
 on $\xi _{2}\mid _{(f\circ \gamma )(J)}$,
i.e. $C(s;\gamma ,f\circ \gamma ):=C_{\gamma (s)}:=
C\mid \pi ^{-1}_{1}(\gamma (s)):\pi ^{-1}_{1}(\gamma (J)) \to
\pi ^{-1}_{2}((f\circ \gamma )(J)), \quad s\in J$;
 b) $(C,f)$ is consistent with $({}^{1}\! I,{}^{2}\! I)$ and c) the equality (2.3$'$)
 is valid.
 	\par {\bf Proof.} If $(F,f)$ and $({}^{1}\! I,{}^{2}\! I)$ are consistent
 (resp. along $\gamma $), then by proposition 2.1 and corollary 2.1 (2.3$'$) and
(2.5) are fulfilled; so defining the bundle morphism $(C,f)$ through
$C_{\gamma (s_{0})}:=C\mid \pi ^{-1}_{1}(\gamma (s_{0}))
:=C(s_{0};\gamma ,f\circ \gamma ), \quad s\in J$ from (2.5) and definition 2.1
(see also (2.4)), we see that $(C,f)$ and $({}^{1}\! I,{}^{2}\! I)$ are consistent
(resp. along $\gamma $). On the opposite, if there exists a bundle morphism
$(C,f)$ with the pointed properties, then there are valid (2.3$'$) (the condition
c)) and (2.5) (follows from the condition b) and (2.4) which is
equivalent to definition 2.1) and by proposition 2.1 $(F,f)$ and
$({}^{1}\! I,{}^{2}\! I)$ are consistent (resp. along $\gamma $).\hbox{\vrule width 6pt height 6pt depth 0pt}
	\par From a functional point of view the general structure
 of the transports
along paths is described by [2], theorem 3.1 and has the form (1.3). The usage
of this theorem allows us to clear up the sense of proposition 2.1, as
well as to solve (locally, i.e. along a given path) the question for
the full description of all bundle morphisms which are (locally) consistent
with a given pair of transports along paths.
	\par Let, in accordance with
[2], theorem 3.1, be chosen sets $Q_{1}$ and $Q_{2}$ and one-to-one maps
${}^{h}\! F^{\gamma _{h}}_{s_{h}}:\pi ^{-1}_{h}(\gamma _{h}(s_{h})) \to Q_{h}$,
 $h=1,2$, which are associated, respectively, with the paths
$\gamma _{h}:J_{h} \to B_{h}, \quad s_{h}\in J_{h}, \quad h=1,2$ and are such
that (cf. (1.3))
\begin{equation}
{}^{h}\! I^{\gamma _{h}}_{s_{h}\to t_{h}}=
\left( {{}^{h}\! F^{\gamma _{h}}_{t_{h}} } \right) ^{-1} \circ
{}^{h}\! F^{\gamma _{h}}_{s_{h}}, \quad s_{h}, \> t_{h}\in J_{h}, \quad h=1,2.
\end{equation}
 	\par
	{\bf Proposition 2.2.} The bundle morphism $(F,f)$ and the pair
$({}^{1}\! I,{}^{2}\! I)$ of transports along paths, which are given through (2.6) by
means of the maps ${}^{1}\! F$ and ${}^{2}\! F$, are consistent (resp. along a path
$\gamma $) iff there exists a map
	\begin{equation}
 C_{0}(\gamma ,f\circ \gamma):Q_{1} \to Q_{2},
	\end{equation}
 such that
	\begin{equation}
F_{\gamma (s)}={\left( {}^{2}\! F^{f\circ \gamma }_{s} \right) }^{-1}\circ
C_{0}(\gamma ,f\circ \gamma )\circ \left( {{}^{1}\! F^{\gamma }_{s} } \right),
	\end{equation}
 or, equivalently, that

\begin{equation}
F_{\gamma (s)}={}^{2}\! I^{f\circ \gamma }_{s_{0} \to s}\circ
C(s_{0};\gamma ,f\circ \gamma )\circ {}^{1}\! I^{\gamma }_{s\to s_{0}},
\end{equation}
where $s_{0}\in J$ is arbitrary and
\begin{equation}
C(s_{0};\gamma ,f\circ \gamma ):=
\left( {{}^{2}\! F^{f\circ \gamma }_{s_{0}} } \right) ^{-1}\circ
C_{0}(\gamma ,f\circ \gamma )\circ
{\left( {}^{1}\! F^{\gamma }_{s_{0}} \right) }
\end{equation}
for every (resp. the given) path $\gamma $.
	\par {\bf Proof.} The substitution of (2.6) into (2.1) shows the
equivalence of the latter with
	\begin{displaymath}
\left( {{}^{2}\! F^{f\circ \gamma }_{t} } \right) \circ
F_{\gamma (t)}\circ \left( {{}^{1}\! F^{\gamma }_{t} } \right) ^{-1}
=\left( {{}^{2}\! F^{f\circ \gamma }_{s} } \right) \circ
F_{\gamma (s)}\circ \left( {{}^{1}\! F^{\gamma }_{s} } \right) ^{-1}
	\end{displaymath}
for any $s,t\in $J.
Hence, if $(F,f)$ and $({}^{1}\! I,{}^{2}\! I)$ are consistent (resp. along
$\gamma $), i.e. (2.1) is satisfied, then the last expression does not at
all depend on $s,t\in J$ and, e.g., fixing arbitrary some $t_{0}\in J$ and
putting $t=t_{0}$ and $C_{0}(\gamma ,f\circ \gamma ):=\left( {{}^{2}\! F^{f\circ
\gamma }_{t_{0}} } \right) \circ F_{\gamma (t_{0})}\circ
\left( {{}^{1}\! F^{\gamma }_{t_{0}} } \right) ^{-1}$ from the last equality, we
easily obtain (2.8).

	\par The equivalence of (2.8) and (2.9) follows directly
from the eq. (2.6) and definition (2.10).
	\par On the opposite, if it is
valid (resp. along $\gamma $), e.g., (2.9), then by proposition 2.1 the
bundle morphism $(F,f)$ and the pair $({}^{1}\! I,{}^{2}\! I)$ of transports along
paths are consistent (resp. along $\gamma $).
\hbox{\vrule width 6pt height 6pt depth 0pt}
	\par Evidently, the
proposition 2.2 is a direct generalization of proposition 2.2 of [1], which
is its special case.
	\par The difference between propositions 2.1 and 2.2
is that the latter, through the equality (2.10), establishes the
general functional form of the map (2.2).
	\par As by [2], proposition 3.5
the maps ${}^{h}\! F^{\gamma _{h}}_{s_{h}}, \quad h=1,2$ are defined up to a
transformation of a form (see [2], eq. (3.11))
${}^{h}\! F^{\gamma _{h}}_{s_{h}} \to \left( {{}^{h}\! D^{\gamma _{h}} }\right) \circ
\left( {{}^{h}\! F^{\gamma _{h}}_{s_{h}} } \right), \quad h=1,2$, where
${}^{h}\! D^{\gamma _{h}}:Q_{h} \to Q^{\circ }_{h}$, is a one-to-one map of
$Q_{h}$ onto some set $Q^{\circ }_{h}, \quad h=1,2,$ then there exists also
nonuniqueness in the choice of the map (2.8). An elementary check shows the
validity of (cf. [1], eq. (2.7))
\begin{eqnarray}
& {}^{h}\! F^{\gamma _{h}}_{s_{h}} \to \left( {{}^{h}\! D^{\gamma _{h}} }\right) \circ
\left( {{}^{h}\! F^{\gamma _{h}}_{s_{h}} } \right), \quad h=1,2 & \nonumber \\
 & \Longleftrightarrow
C_{0}(\gamma ,f\circ \gamma ) \to \left( {{}^{2}\! D^{f\circ \gamma } }
\right) ^{-1}\circ C_{0}(\gamma ,f\circ \gamma )\circ
 \left( {{}^{1}\! D^{\gamma } } \right). &
\end{eqnarray}

	\par {\bf Proposition 2.3.} If for a given pair
$({}^{1}\! I,{}^{2}\! I)$ of transports along
paths the representation (2.6) is chosen, then all consistent along
$\gamma :J\to B_{1}$ with it bundle morphisms $(F^{\gamma ,f\circ \gamma },f)$
along $\gamma $ are obtained from the equality
\begin{equation}
F^{\gamma ,f\circ \gamma }_{\gamma (s)}:=
 {\left( {{}^{2}\! F^{f\circ \gamma }_{s} } \right) }^{-1}\circ
 C_{0}(\gamma ,f\circ \gamma )\circ \left( { {}^{1}\! F^{\gamma }_{s} }
 \right) =  {}^{2}\! I^{f\circ \gamma }_{s_{0} \to s}\circ
 C(s_{0};\gamma ,f\circ \gamma )\circ {}^{1}\! I^{\gamma }_{s \to s_{0}},\
\end{equation}
where $s_{0}\in J$ is arbitrary, $C$ is
defined by (2.10) in which $C_{0}(\gamma ,f\circ \gamma )$ is a
one-to-one map from $Q_{1}$ onto $Q_{2}$.
	\par {\bf Proof.} This proposition is a
consequence of the proof of proposition 2.2, as from it is clear that
(2.12) is the general solution of the equation (2.1) with respect to
$F_{\gamma (t)}$ when ${}^{1}\! I$ and ${}^{2}\! I$ are given.
\hbox{\vrule width 6pt height 6pt depth 0pt}
	\par The
definition 1.1 of [1] for consistency between bundle metrics and transports
along paths seems rather natural by itself for a difference of the
definition 2.1 for consistency between bundle morphisms and a pair of
transports along paths, whose introduction needs some explanation. As we
saw in the written after definition 2.1 the former definition is a
special case of the latter. Now we shall show that in this context the
definition 2.1 itself is a special case of definition 2.2 of [2]
for a section of a given fibre bundle transported along a path.
	\par
	Let there be given two fibre bundles
$\xi _{h}=(E_{h},\pi _{h},B_{h}), \quad h=1,2.$ We define the fibre bundle
$\xi _{0}=(E_{0},\pi _{0},B_{1})$ of bundle morphisms from $\xi _{1}$
on $\xi _{2}$ in the following way:
\begin{eqnarray}
&\!\! \! \! \! \! \! \! \! \! \!
E_{0}:=\{ (F_{b_{1}},f)\! : F_{b_{1}}\! \! :\pi ^{-1}_{1}(b_{1}) \to
	\pi ^{-1}_{2}(f(b_{1})), \  b_{1}\in B_{1}, \
		f:B_{1} \to B_{2} \}, \quad   & \\
 & \pi _{0}((F_{b_{1}},f)):=b_{1}, \quad (F_{b_{1}},f)\in E_{0}, \quad
 b_{1}\in B_{1}. &
\end{eqnarray}
	\par It is clear that every section $(F,f)\in Sec\xi _{0}$
is a bundle morphism from $\xi _{1}$ into $\xi _{2}$ and vice versa,
every bundle morphism from $\xi _{1}$
on $\xi _{2}$ is a section of $\xi _{0}$. (Thus a bundle structure in the
set $Morf(\xi _{1},\xi _{2})$ of bundle morphisms from
 $\xi _{1}$ on $\xi _{2}$ is introduced.)
	\par
	If in $\xi _{0}$ given is a transport $K$ along the paths
in $B_{1}$, then
according to [2], definition 2.2 (see therein eq. (2.4)), the bundle
morphism $(F,f)\in Sec\xi _{0}$ is $(K-)$transported along
 $\gamma :J \to B_{1}$ if
\begin{equation}
 (F_{\gamma (t)},f)=K^{\gamma }_{s \to t}(F_{\gamma (s)},f), \quad s,t\in J.
\end{equation}
	\par
	If in $\xi _{1}$ and $\xi _{2}$ are given, respectively,
the transports ${}^{1}\! I$ and ${}^{2}\! I$ along the paths, respectively,
in $B_{1}$ and $B_{2}$, then they generate
in $\xi _{0} $ a "natural" transport $^{0}\! K$ along the paths in $B_{1}$.
The action of this transport along
$\gamma :J \to B_{1}$ on $(F_{\gamma (s)},f)\in \pi ^{-1}_{0}(\gamma (s))$
for a fixed $s\in J$ and arbitrary $t\in J$ is defined by
	\begin{equation}
 {}^{0}\! K^{\gamma }_{s \to t}(F_{\gamma (s)},f):=
 \left( {{}^{2}\! I^{f\circ \gamma }_{s \to t} \circ
 F_{\gamma (s)}\circ {}^{1}\! I^{\gamma }_{t \to s}, \; f} \right) \in
 \pi ^{-1}_{0}(\gamma (t)).
	\end{equation}
	\par
	 {\bf Proposition 2.4.} The  map
${}^{0}\! K^{\gamma }_{s \to t}:\pi ^{-1}_{0}(\gamma (s)) \to
 \pi ^{-1}_{0}(\gamma (t))$
defined through (2.16)  is a transport along $\gamma $ from $s$ to $t,$
 $s,t\in J$ and, consequently, in $\xi _{0} \quad ^{0}\! K$ is a transport
along paths.
	\par
	{\bf Proof.} Using the properties (1.1) and (1.2) of the transports
along paths it is easy to check with the help of (2.16) that the
maps $^{0}\! K^{\gamma }_{s \to t}, \quad s,t\in J$ satisfy the equalities
		\begin{eqnarray}
 & {}^{0}\! K^{\gamma }_{s \to t}\circ {}^{0}\! K^{\gamma }_{r \to s}=
	{}^{0}\! K^{\gamma }_{r\to t}, \quad r,s,t\in J, & \\
 & {}^{0}\! K^{\gamma }_{s \to s}={\it id}_{\pi ^{-1}_{0}(\gamma (s))},
\quad s\in J &
		\end{eqnarray}
and hence, by [2], definition 2.1, ${}^{0}\! K^{\gamma }_{s \to t}$ is a
transport along $\gamma $ from $s$ to $t$, i.e. in $\xi _{0} \  {{}^{0}\! K}$
really defines a transport along paths.\hbox{\vrule width 6pt height 6pt
	depth  0pt}
	\par
	{\bf Lemma 2.1.} If $(F,f)\in Sec\xi _{0}$, then (2.1) is  equivalent to
		\begin{equation}
	(F_{\gamma (t)},f)={}^{0}\! K^{\gamma }_{s \to
t}(F_{\gamma (s)},f), \quad s,t\in J.
		\end{equation}
		\par
		{\bf Proof.} At the
	beginning of the proof of proposition 2.1 we saw that (2.1) is
 equivalent to (2.4), that  is equivalent to (2.19) because of (2.16), i.e.
(2.1) and (2.19) are equivalent.\hbox{\vrule width 6pt height 6pt depth 0pt}
	\par {\bf Proposition 2.5.} The bundle morphism $(F,f)$ and
the pair $({}^{1}\! I,{}^{2}\! I)$ of transports along paths are
consistent (resp. along the
path $\gamma $) iff $(F,f)$ is transported along every (resp. the given) path
$\gamma $ with the help of the defined from $({}^{1}\! I,{}^{2}\! I)$
in $\xi _{0}$ transport along paths $^{0}\! K$ .
	\par {\bf Proof.} The proposition follows directly from
lemma 2.1, definition 2.1 and definition 2.2 of [2] (see (2.15) and
[2], eq. (2.4)).\hbox{\vrule width 6pt height 6pt depth 0pt}
	\par Taking into account the comment after definition 2.1 it is not
difficult to verify that proposition 3.1$'$ of [1] is, in fact,
a variant of proposition 2.5 in the special case when one studies
the consistency of $S$-transports and bundle metrics.

\setcounter {chapter} {3}
\setcounter {equation} {0}
\section {\bf CONSISTENCY WITH MORPHISMS \\ OF THE FIBRE BUNDLE}

	In this section we are going to apply the general theory of
Sect.~2 to the case of bundle morphisms of a given fibre bundle.
The results so-obtained will be illustrated with often used examples.
\par
 Let in the fibre bundle $(E,\pi ,B)$ a (bundle) morphism $F$ be given.
By definition [4,5] this means that $(F,{\it id}_{B})$ is a bundle morphism
from $(E,\pi ,B)$ into $(E,\pi ,B)$. Hence it is natural (cf. definition 2.1)
$F$ and the transport along paths $I$ in $(E,\pi ,B)$ to be called globally
(resp. locally) consistent if $(F,{\it id}_{B})$ and the pair $(I,I)$ are
 globally (resp. locally) consistent, i.e. in this case definition 2.1
reduces to
	\par
	{\bf Definition 3.1.} The transport along paths $I$ in $(E,\pi ,B)$
is globally (resp. locally) consistent with the bundle morphism $F$ of
$(E,\pi ,B)$ if
\begin{equation}
F_{\gamma (t)}\circ I^{\gamma }_{s\to t}=I^{\gamma }_{s\to t}\circ
 F_{\gamma (s)}, \quad s,\, t\in J
\end{equation}
is fulfilled for every (resp. the given) path $\gamma :J \to B$.
\par
This definition formalizes the condition for commutation of a
transport along paths and a bundle morphism of the fibre bundle, in
which the transport acts, or, in other words, the equality
(3.1) is an exact expression of the phrase that "the bundle morphism
$F$ and the transport along paths $I$ commute".
	\par
	Comparing definitions 3.1 and 2.1, we see that from the first
of them  the second can be obtained if in the former we put
	\begin{equation}
 \left( {E_{h},\pi _{h},B_{h} }\right) =(E,\pi ,B), \quad h=1,2,
 \  f={\it id}_{B}, \ {}^{1}\! I={}^{2}\! I=I.
	\end{equation}

	If we make these substitutions in the whole section 2 and take
into account definition 3.1, then the stated therein propositions
and definitions, concerning bundle morphisms between fibre bundles,
take the following formulations in the case of bundle morphisms of
the fibre bundle $(E,\pi ,B)$:
\par
	{\bf Proposition 3.1.} The bundle morphism $F$ and the transport along
paths $I$ are consistent (resp. along $\gamma $) if and only if there exist
$s_{0}\in J$ and a map
\begin{equation}
C(s_{0};\gamma ):\pi ^{-1}(\gamma (s_{0})) \to \pi ^{-1}(\gamma (s_{0})),
\end{equation}
such that
\begin{equation}
F_{\gamma (s)}=I^{\gamma }_{s_{0}\to s}\circ
C(s_{0};\gamma )\circ I^{\gamma }_{s\to s_{0}}
\end{equation}
for every (resp. the given) path $\gamma $.
	\par
{\bf Corollary 3.1.} The equality (3.4) is true iff for every $t_{0}\in J$,
we have
		\setcounter{equation}{3}
		\def\theequation{\thechapter.\arabic{equation}'}
\begin{equation}
F_{\gamma (s)}=I^{\gamma }_{t_{0}\to s}\circ
C(t_{0};\gamma )\circ I^{\gamma }_{s\to t_{0}},
\end{equation}
where
		\setcounter{equation}{4}
		\def\theequation{\thechapter.\arabic{equation}}
\begin{equation}
C(t_{0};\gamma )=I^{\gamma }_{s_{0}\to t_{0}}\circ
 C(s_{0};\gamma )\circ I^{\gamma }_{t_{0}\to s_{0}}.
\end{equation}
	\par
	{\bf Corollary 3.2.} The bundle morphism $F$ and the transport along
paths $I$ are consistent (resp. along $\gamma $) iff there exists a bundle
morphism $C$ of $(E,\pi ,B)$ such that for every (resp. the given) path
$\gamma :J\to B$ we have:
a) $C\mid \pi ^{-1}(\gamma (J))$ is a morphism of
 $(E,\pi ,B)\mid _{\gamma (J)}
=\left( {\pi ^{-1}(\gamma (J)),\pi \mid \gamma (J),\gamma (J) }\right) $;
 b) $C$ and $I$ are consistent (resp. along $\gamma $);
 and c) (3.4$'$) is true.
	\par
	{\bf Proposition 3.2.} The bundle morphism $F$ and
 the transport  along paths $I$ defined by the equality
$I^{\gamma }_{s \to t}=\left( F^{\gamma }_{t} \right) ^{-1}\circ
F^{\gamma }_{s}, \quad t,\, s\in J $ (see (1.3)) are consistent (resp. along $\gamma $) if and only if there exists a map
$C_{0}(\gamma ):Q \to Q$, such that
\begin{equation}
F_{\gamma (s)}=\left( F^{\gamma }_{s} \right) ^{-1}\circ
C_{0}(\gamma )\circ F^{\gamma }_{s}=I^{\gamma }_{t \to s}\circ
C(s;\gamma )\circ I^{\gamma }_{s \to t}, \quad s,t\in J
\end{equation}
where
\begin{equation}
C(s;\gamma ):=\left( F^{\gamma }_{s} \right) ^{-1}\circ
 C_{0}(\gamma )\circ F^{\gamma }_{s}.
\end{equation}
	\par
	{\bf Proposition 3.3.} If a transport along paths $I$ with a
representation
$I^{\gamma }_{s \to t}=\left( F^{\gamma }_{t}\right) ^{-1}\circ
F^{\gamma }_{s}, \quad s\in J $ (see (1.3)) is fixed, then any
bundle morphism $^{\gamma }\! F$ along $\gamma $  consistent along
 $\gamma :J \to B$ with it is obtained from the equality
\begin{equation}
 ^{\gamma }\! F_{\gamma (s)}:=\left( F^{\gamma }_{s}\right) ^{-1}\circ
C_{0}(\gamma )\circ F^{\gamma }_{s}=I^{\gamma }_{s_{0} \to s}\circ
C(s_{0};\gamma )\circ I^{\gamma }_{s \to s_{0}},
\end{equation}
where $s_{0}\in J$ is arbitrary, $C_{0}:Q \to Q$ is one-to-one and $C$ is
defined by (3.7).
\par
According to (2.13) and (2.15) the fibre bundle $(E_{0},\pi _{0},B)$ of
the bundle morphisms of $(E,\pi ,B)$ is defined through the equalities
\begin{eqnarray}
& E_{0}:=\left\{  F_{b}: \quad F_{b}:\pi ^{-1}(b) \to \pi ^{-1}(b),
 \quad b\in B \right\} , & \\
& \pi _{0}(F_{b}):=b, \quad F_{b}\in E_{0}, \quad b\in B. &
\end{eqnarray}
Evidently, if $F$ is a morphism of $(E,\pi ,B)$, then
$F\in Sec(E_{0},\pi _{0},B)$ and vice versa.
\par
According to [2], definition 2.2 if in $(E_{0},\pi _{0},B)$ is given a
transport along paths $K$, then $F\in Sec(E_{0},\pi _{0},B)$ is
$K$-transported along $\gamma :J \to B$ if
\begin{equation}
F_{\gamma (t)}=K^{\gamma }_{s \to t}(F_{\gamma (s)}).
\end{equation}

	If $I$ is a transport along paths in $(E,\pi ,B)$, then in
$(E_{0},\pi _{0},B)$ it induces in
 a "natural" transport along paths $^{0}\! K$ whose action
along
$\gamma :J \to B$ on $F_{\gamma (s)}\in \pi ^{-1}_{0}(\gamma (s))$ is
	\begin{equation}
^{0}\! K^{\gamma }_{s \to t}(F_{\gamma (s)}):=I^{\gamma }_{s \to t}\circ
F_{\gamma (s)}\circ I^{\gamma }_{t \to s}\in \pi ^{-1}_{0}(\gamma (t)),
\quad  s,t\in J.
	\end{equation}
	\par
	{\bf Proposition 3.4.} The mapping
$^{0}\! K^{\gamma }_{s \to t}:\pi ^{-1}_{0}(\gamma (s)) \to
\pi ^{-1}_{0}(\gamma (t))$,
defined by (3.12), is a transport along $\gamma $ from $s$ to $t, \quad
s,t\in J$ and, consequently, $^{0}\! K$ defines a
 transport along paths in $(E_{0},\pi _{0},B)$ .
	\par
	{\bf Lemma 3.1.} The equality (3.1) is equivalent to
	\begin{equation}
F_{\gamma (t)}={}^{0}\! K^{\gamma }_{s \to t}(F_{\gamma (s)}), \quad s,t\in J.
	\end{equation}
	\par
	{\bf Proposition 3.5.} In $(E_{0},\pi _{0},B)$ the bundle morphism
$F$ and the
transport along paths $I$ are consistent (resp. along $\gamma $) iff $F$ is
transported along every (resp.  the given) path $\gamma $ with the help of
the transport along paths $^{0}\! K$ defined by $I$ in $(E_{0},\pi _{0},B)$.
\par
Now we shall consider examples for consistency of concrete bundle
 morphisms in vector bundles with transport along paths in them.
	\par
	\vspace{3mm} {\bf\Large Example 3.1.} {\large\sl Consistency with an almost complex structure}
	\par
Let the bundle morphism ${\bf J}$ of the real vector bundle $(E,\pi ,B)$
define an almost complex structure in it [8-10], i.e.
${\bf J}_{x}:={\bf J}\mid {\pi }^{-1}(x), \quad x\in B$ to be
${\bf R}$-linear
isomorphisms defining complex structure in the fibres $\pi ^{-1}(x)$, which
means that
\begin{equation}
 {\bf J}_{x}\circ {\bf J}_{x}:=-{\it id}_{\pi ^{-1}(x)}.
\end{equation}

Evidently, if $(E,\pi ,B)$ is the tangent bundle to some manifold and
${\bf J}$ is a linear endomorphism, then ${\bf J}$ defines
an almost complex structure on that manifold [3,8].
In this case, following the accepted
terminology, a transport along paths consistent with ${\bf J}$ may be called
almost complex.
	\par
	{\bf Proposition 3.6.} A bundle morphism ${\bf J}$ consistent with a
transport along paths $I$ of the vector fibre bundle
$(E,\pi ,B)$ defines an almost complex structure in it if and only if the
involved in (3.4)-(3.7) (with ${\bf J}_{\gamma (s)}$ instead of
$F_{\gamma (s)}$) map $C(s;\gamma )$ or map $C_{0}(\gamma )$ define a complex
structure in $\pi ^{-1}(\gamma (s))$ or $Q$ respectively, i.e. when they
satisfy the following equalities:
	\begin{equation}
 C(s;\gamma )\circ C(s;\gamma )=-{\it id}_{\pi ^{-1}(\gamma (s))},
	\end{equation}
		\setcounter{equation}{14}
		\def\theequation{\thechapter.\arabic{equation}'}
	\begin{equation}
 C_{0}(\gamma )\circ C_{0}(\gamma )=-{\it id}_{Q}.
	\end{equation}
		\setcounter{equation}{15}
		\def\theequation{\thechapter.\arabic{equation}}Or,
 in other words, the almost complex structure ${\bf J}$
and the transport along paths $I$ are globally (resp.  locally)
consistent, i.e. they commute globally (resp.  locally), if and only if
there are fulfilled (3.4) (or (3.6)), (3.15) and (3.15$'$) for every (resp.
the given) path $\gamma $; besides, eqs.  (3.15) and (3.15$'$) are equivalent.
	\par
	{\bf Proof.} According to proposition 3.1 ${\bf J}$
and $I$ are consistent (resp.  along $\gamma $)
iff (3.4) holds for $F={\bf J}$ and every (resp. the given) path
$\gamma $, so we get
${\bf J}_{\gamma (s)}\circ {\bf J}_{\gamma (s)}=
I^{\gamma }_{s_{0} \to s}\circ C(s_{0};\gamma )\circ C(s_{0};\gamma )\circ
 I^{\gamma }_{s \to s_{0}}$,
whence it follows that (resp.  along $\gamma $) (3.13) is
equivalent to (3.15) for $s=s_{0}$.
	\par Analogously, based on the
considerations on (3.6), one proves (resp. along $\gamma $) the equivalence of
(3.13) with (3.14) (for every $s$) and (3.15$'$).
  \par The equivalence
of (3.15) and (3.15$'$) is a consequence of
(3.7).\hbox{\vrule width 6pt height 6pt depth 0pt}
	\par
\vspace{3mm} {\bf\Large Example 3.2.} {\large\sl Consistency with a multiplication with numbers}
\par Let $(E,\pi ,B)$ be a real (resp. complex) vector bundle,
$\lambda \in {\bf R}$ (resp.  $\lambda \in {\bf C}$) and the bundle morphism
$^{\lambda }\! F$ of $(E,\pi ,B)$ be defined by
 $^{\lambda }\! F(u):=\lambda \cdot u={ } ^{\lambda }\! F_{\pi (u)}(u)$ for every
$u\in $E.
\par
	{\bf Definition 3.2.} The transport along paths $I$ in the real
(resp.  complex) vector bundle
$(E,\pi ,B)$ is called consistent with the operation multiplication with
real (resp.  complex) numbers if it is consistent with the bundle morphisms
$^{\lambda }\! F$ for every $\lambda \in {\bf R}$ (resp. $\lambda \in {\bf C})$.
	\par {\bf Proposition 3.7.} The transport along paths $I$
is globally (resp.
locally) consistent with the multiplication with, respectively, real or
complex numbers if and only if
\begin{equation}
I^{\gamma }_{s \to t}(\lambda u)=\lambda
\left( I^{\gamma }_{s \to t}(u)\right) , \quad u\in \pi ^{-1}(\gamma (s))
\end{equation}
for every, respectively, $\lambda \in {\bf R}$
or $\lambda \in {\bf C}$ and every (resp. the given) path $\gamma :J \to B$.
\par
{\bf Proof.} The proposition is a simple corollary of definitions 3.1
and 3.2, proposition 3.5, lemma 3.1
and (3.12).\hbox{\vrule width 6pt height 6pt depth 0pt}
	\par
In other words we may say that the consistency with multiplication
with numbers means simply the validity of the condition for homogeneity (3.16)
or, which is the same, the operations of $I$-transportation along paths and
multiplication with numbers in the fibres to commute.
\par
\vspace{3mm} {\bf\Large Example 3.3.} {\large\sl Consistency with the operation addition}
\par
Let $(E,\pi ,B)$ be a vector fibre bundle, $A\in Sec(E,\pi ,B)$
and $^{A}\! F$
be a bundle morphism of $(E,\pi ,B)$ defined by
$^{A}\! F(u):=u+A(\pi (u))=^{A}\! F_{\pi (u)}(u)$
for every $u\in E$.
 	\par
	{\bf Definition 3.3.} The
transport along paths $I$ in the vector fibre bundle $(E,\pi ,B)$ is called
globally (resp.  locally) consistent with the operation addition if it is
consistent with the bundle morphisms $^{A}\! F$ for every section
A $I$-transported along every (resp.  the given) path.
	\par
	{\bf Proposition 3.8.} The transport $I$ is globally
(resp.  locally) consistent with the operation
addition iff these operations commute, i.e. iff
\begin{equation}
I^{\gamma }_{s \to t}(u+v)=I^{\gamma }_{s \to t}(u)+I^{\gamma }_{s \to t}(v),
\quad u,v\in \pi ^{-1}(\gamma (s)),
\end{equation}
 which means that the usual additivity
condition in the fibres is to be fulfilled for every (resp.  the given) path
$\gamma :J \to $B.
	\par
	{\bf Proof.} According to definition 3.1 the
transport $I$ and the above defined bundle morphisms $^{A}\! F$ are globally
(resp. locally) consistent iff the equality
\begin{equation}
A(\gamma (t))+I^{\gamma }_{s \to t}(u)=I^{\gamma }_{s \to t}[A(\gamma (s))+u],
\quad  u\in \pi ^{-1}(\gamma (s))
\end{equation}
is valid for every (resp.  the given) path
$\gamma $.
	\par
	The condition A to be I-transported along $\gamma $ section
means (see [2], definition 2.2 and proposition 2.1) that $A(\gamma (t))=
I^{\gamma }_{s \to t}A(\gamma (s)), \quad s,t\in J$, which, when substituted
into the previous equality, due to the arbitrariness of A (i.e. of $A(\gamma
(s))$) gives (3.17) (see [2], proposition 2.2) in which $A(\gamma (s))$ is
denoted with v.  On the opposite, if (3.17) is valid and A is a section
I-transported along $\gamma $ section, then (3.18)
is identically satisfied and, consequently, $^{A}\! F$ and $I$
are consistent.\hbox{\vrule width 6pt height 6pt depth 0pt}
 	\par
\vspace{3mm} {\bf\Large Example 3.4.} {\large\sl Consistency between transport along paths and Finslerian metrics}
	\par
Let in a manifold $M$ a Finslerian metric be given [11] by means of
a Finslerian metric function
$F:T(M) \to {\bf R}_{0}:=\{ \lambda : \ \lambda \in {\bf R}, \
\lambda \ge 0\} $ having the property $F(x,\lambda A)=\lambda F(x,A)$ for
$\lambda \in {\bf R}_{+}:=\{ \lambda : \  \lambda \in {\bf R},\
\lambda >0\} $, $x\in M, \ A\in T_{x}(M)$
and satisfying the conditions described in the above references.  Here
$T(M):=\cup_{x\in M} T_{x}(M)$, where $T_{x}(M)$
is the tangent to $M$ space at $x\in $M.
	\par {\bf Definition 3.4.} The
Finslerian metric and the transport along paths $I$ are consistent (resp.
along $\gamma :J \to M$) if the equality
\begin{equation}
F(\gamma (s),A)=F(\gamma (t),I^{\gamma }_{s \to t}A), \quad s,t\in J,
\quad A\in T_{\gamma (s)}(M)
\end{equation}
is fulfilled for every (resp.  the given) path $\gamma $.
\par This definition is a special case of definition
2.1 and it is obtained from it for:
$\xi _{1}=(T(M),\pi ,M)$; $\xi _{2}=({\bf R}_{0},\pi _{0},0)$, where
$0\in {\bf R}$ and $\pi _{0}:{\bf R}_{0} \to \{ 0\} $; in the bundle morphism
$(F,f) \quad F:T(M) \to {\bf R}_{0}$ is the Finslerian metric function,
$F_{x}(A)=F(x,A), \quad A\in T_{x}(M)$ and $f:M \to \{0\}$;
 ${}^{1}\! I^{\gamma }=I^{\gamma }$;
${}^{2}\! I^{f\circ \gamma }={\it id}_{{\bf R}_{0}}$.
 \par \vspace{3mm} {\bf\Large Example 3.5.} {\large\sl Consistency between transports along paths and
symplectic metrics}
\par A real symplectic metric $a:x\mapsto a_{x}, \quad
a_{x}:\pi ^{-1}(x)\times \pi ^{-1}(x) \to {\bf R}$ in a fibre
bundle $(E,\pi ,B)$
differs from a real symmetric metric $g$ (cf.  [1]) only in that it is
antisymmetric, i.e.  $a_{x}(u,v)=-a_{x}(v,u)$ for $u,v\in \pi ^{-1}(x)$.
Hence, modifying definition 1.1 of [1], we can say that the transport along
paths $I$ and $a$ are consistent (resp.  along $\gamma $) if
	\begin{equation}
a_{\gamma (s)}=a_{\gamma (t)}\circ (I^{\gamma }_{s \to t}\times
I^{\gamma }_{s \to t}), \quad s,t\in J
	\end{equation}
is fulfilled for every (resp.  the given) path $\gamma :J \to $B.

\setcounter {chapter} {4}
\setcounter {equation} {0}
\section {\bf  CONSISTENCY WITH A HERMITIAN \\ STRUCTURE}

 By a Hermitian structure in the real vector bundle
$(E,\pi ,B)$ we understand (cf. [10]) a pair $({\bf J},g)$ of almost complex
structure ${\bf J}$, i.e.  a bundle morphism ${\bf J}:E \to E$ with the
property ${\bf J\circ J}$=-{\it id}$_{E}$, and a consistent with it bundle
symmetric metric $g \  (g:x \to g_{x}$, such that $g_{x}$ is bilinear,
symmetric and $g_{x}=g_{x}\circ ({\bf J}_{x}\times {\bf J}_{x})$) in that
bundle (called a Hermitian metric; see [10], ch.  IX, \S 1).
	\par
	{\bf Definition 4.1.} If $I$ is a transport along paths in
$(E,\pi ,B)$, then the
Hermitian structure $({\bf J},g)$ and $I$ are consistent (resp.  along the
path $\gamma $) if the pairs ${\bf J}$ and $I$ and $g$ and $I$ are consistent
separately, i.e.
	\begin{eqnarray}
& {\bf J}_{\gamma (t)}\circ I^{\gamma }_{s \to t}=
I^{\gamma }_{s \to t}\circ {\bf J}_{\gamma (s)}, &\\
& g_{\gamma (s)}=g_{\gamma (t)}\circ
\left( I^{\gamma }_{s \to t}\times I^{\gamma }_{s \to t} \right), &
	\end{eqnarray}
for every (resp.  the given) path $\gamma $.
	\par
	{\bf Remark.} The condition
(4.2) for consistency between $g$ and $I$ was introduced and investigated in
[1].
\par In particular, if $(E,\pi ,B)=(T(M),\pi ,M))$ is the tangent bundle
to the manifold $M$, then (4.1) and (4.2) define (an "almost Hermitian")
transport consistent with the Hermitian structure $({\bf J},g)$ of the almost
Hermitian manifold $(M,{\bf J},g)$ (cf.  [10], ch.  IX, \S 2).
 	\par
	Further we shall dwell on the question for consistency between linear
transports ($L$-transports) along paths $L$ in a vector bundle $(E,\pi ,B)$
[7] and Hermitian structures $({\bf J},g)$ in it, such that ${\bf J}$ is a
linear endomorphism.
	\par
	Let along the path $\gamma :J \to B$ there be fixed bases
$\{e_{i}(s), \quad i=1,\ldots ,$ $\dim (\pi ^{-1}(x))$, $x\in B\}$ in
$\pi ^{-1}(\gamma (s)), \quad s\in J$, in which ${\bf J}$ is defined by the
matrices ${\bf J}(\gamma (s))= \| {\bf J}^{i}_{.j}(\gamma (s))\| $ , $L$
is defined through the matrices $F(s;\gamma )$ by (1.4) and $g$ - through the
matrix $G(\gamma (s))=\Vert g_{\gamma (s)}(e_{i}(s),e_{j}(s)) \Vert $. Then
the condition $g=g\circ ({\bf J\times J})$ for consistency along $\gamma $
between ${\bf J}$ and $g$ takes the form
\begin{equation}
 {\bf J}^{\top }(\gamma (s))G(\gamma (s)){\bf J}(\gamma (s))=G(\gamma (s)),
\end{equation}
where $\top $ means transposition of matrices, and the the condition
for consistency between
${\bf J}$ and $I$, due to propositions 3.2 and 3.6, looks like
	\begin{equation}
{\bf J}(\gamma (s))=F^{-1}(s;\gamma )C_{0}(\gamma )F(s;\gamma ),
\qquad C_{0}(\gamma )C_{0}(\gamma )=-{\bf I},
	\end{equation}
 where $C_{0}(\gamma )$ is a nondegenerate matrix and $\bf I$ is the unit
 matrix.
	\par Let in
$(E,\pi ,B)$ an $L$-transport along paths $L$ be given [7]. The following two
propositions solve the problems for the existence and a full (local and
global) description of all consistent with $L$ Hermitian structures.
	\par
	{\bf Proposition 4.1.} Along the path $\gamma :J \to B$
the class of all Hermitian structures (along $\gamma $) which are
locally consistent with the
$L$-transport along paths $L$  is given in the above pointed bases
through the equalities
		\setcounter{equation}{4}
		\def\theequation{\thechapter.\arabic{equation}a}
	\begin{equation}
 {\bf J}(\gamma (s);\gamma )=F^{-1}(s;\gamma )C_{0}(\gamma )F(s;\gamma ),
	\end{equation}
		\setcounter{equation}{4}
		\def\theequation{\thechapter.\arabic{equation}b}
	\begin{equation}
 G(\gamma (s);\gamma )=F^{\top }(s;\gamma )C(\gamma )F(s;\gamma )
	\end{equation}
		\setcounter{equation}{5}
		\def\theequation{\thechapter.\arabic{equation}}
in which the matrix
functions $C_{0}$ and $C$ satisfy the equations
	\begin{eqnarray}
{\bf I}=C^{\top }C^{-1}=-C_{0}C_{0}=C^{\top }_{0}CC_{0}C^{-1}.
	\end{eqnarray}\
	\par
	{\bf Proof.} The
equalities (4.5a) and (4.5b), and also the first two from (4.6), follow,
respectively, from (4.4) and [1], proposition 2.3 (see therein eq.  (2.8)).
The last equality from (4.6) is obtained by the substitution of (4.5) into
the condition for consistency between ${\bf J}$ and $g$ expressed now by
(4.3).\hbox{\vrule width 6pt height 6pt depth 0pt}
	\par {\bf Proposition 4.2.} Let there be given an
$L$-transport along paths $L$ defined along $\gamma :J \to B$ through (1.4)
by the matrices
	\begin{equation}
F(s;\gamma )=Y(\gamma )Z(s;\gamma )D^{-1}(\gamma (s)), \quad s\in J,
	\end{equation}
where $Y$ and $D$ are nondegenerate matrix
functions and $Z$ is a pseudoortho\-gonal of some type $(p,q), \
p+q=\dim (\pi ^{-1}(\gamma (s)))$ matrix function, i.e.
	\begin{equation}
Z^{\top }(s;\gamma )G_{p,q}Z(s;\gamma )=G_{p,q}:=
 diag(\underbrace {1,\ldots ,1}_{\small {p-times}},
\underbrace {-1,\ldots  ,-1}_{\small {q-times}}).
 	\end{equation}
Then every
Hermitian structure $({\bf J},g)$ globally consistent with $L$, if any, is
given by the equalities
		\setcounter{equation}{8}
		\def\theequation{\thechapter.\arabic{equation}a}
	\begin{equation}
{\bf J}(\gamma (s))=D(\gamma (s))Z^{-1}(s;\gamma )PZ(s;\gamma )
D^{-1}(\gamma (s)),
	\end{equation}
		\setcounter{equation}{8}
		\def\theequation{\thechapter.\arabic{equation}b}
	\begin{equation}
 G(\gamma (s))=(D^{-1}(\gamma (s)))^{\top }G_{p,q}D^{-1}(\gamma (s)),
	\end{equation}
		\setcounter{equation}{9}
		\def\theequation{\thechapter.\arabic{equation}}
where $P$ is a constant with respect to $s$ and $\gamma $ matrix,
which may depend on $p$ and $q$ and is explicitly constructed below (see
(4.11)), and $Z$ besides (4.8) satisfies the condition that $Z^{-1}(s;\gamma
)PZ(s;\gamma )$ depends only on $\gamma (s)$, but not on $s$ and $\gamma $
separately.
	\par
	{\bf Remark.} The necessity of the representation (4.7) for
$F(s;\gamma )$ is a consequence of that we want a globally consistent
with $L$ metric $g$ to exist (see [1], proposition 2.6).
	\par
	{\bf Proof.} The
equality (4.9b) is a corollary from proposition 2.6 from [1]. The
consistency of $L$ and ${\bf J}$ is equivalent to (4.4), due to which
substituting (4.7) into (4.4), we get (4.9b) with
	\begin{equation}
P=Y^{-1}(\gamma )C_{0}(\gamma )Y(\gamma ), \quad C_{0}(\gamma )C_{0}(\gamma )
=-{\bf I}.
	\end{equation}
 As ${\bf J}$ and $g$ must form a Hermitian structure, they
have to be consistent, i.e. (4.3) must be true which, as we shall now prove,
is a consequence of the independence of $P$ of $\gamma $. In fact,
substituting (4.9) into (4.3) and using (4.8), we find
		\setcounter{equation}{10}
		\def\theequation{\thechapter.\arabic{equation}a}
	\begin{equation}
P^{\top }G_{p,q}P=G_{p,q},
	\end{equation}
i.e. $P$ is a pseudoorthogonal
matrix of type $(p,q)$ which, as a result of (4.10), satisfies
		\setcounter{equation}{10}
		\def\theequation{\thechapter.\arabic{equation}b}
	\begin{equation}
PP=-{\bf I}.
	\end{equation}
		\setcounter{equation}{11}
		\def\theequation{\thechapter.\arabic{equation}}

	If we consider (4.11) as a system of equations
with respect to $P$, then, its solution, if any (see below),
is independent of any parameters as it is a function only of ${\bf I}$ and
$G_{p,q}$, i.e. the elements of $P$ are independent of $\gamma $ numbers.
This conclusion is a consequence of the observation that (4.11a) does not
change when it is transposed, i.e. $P^{\top }G_{p,q}P$ is a symmetric, and
(4.11b), due to (4.11a), is equivalent to
 $(G_{p,q}P)^{\top }\equiv P^{\top }G_{p,q}=-G_{p,q}P$, i.e. $G_{p,q}P$ is
 antisymmetric.  Hence
(4.11a) and (4.11b) contain respectively $n(n+1)/2$ and $n(n-1)/2, \quad
n=\dim (\pi ^{-1}(x)), \quad x\in P$, or commonly $n^{2}$, a number of
independent scalar equations for the $n^{2}$ components of P. So, if $P$
exists, it is constant (along $\gamma $).
  	\par The condition for that the matrix
$Z^{-1}(s;\gamma )PZ(s;\gamma )$ to depend only on $\gamma (s)$ is a result
from that ${\bf J}$ must be globally defined, i.e. ${\bf J}(\gamma (s))$ must
depend only on the point $\gamma (s)$, but not on the path $\gamma :J \to B$,
so from (4.9a) the above pointed
condition follows.\hbox{\vrule width 6pt height 6pt depth 0pt}
	\par {\bf Remark.} It may be proved that for an even $n$, i.e.
for $n=2k+1$, $k=1,2,\ldots $, the equations (4.11) have no solutions with
 respect to $P$, which is in
accordance with the fact that in this case in $(E,\pi ,B)$
a complex structure cannot be introduced
(see [10], ch.  IV, \S 1).  For an odd $n$, i.e.
for $n=2k$, the equations (4.11) have different solutions with respect to P.
For instance, for $n=2$ these solutions for $p=2=2-q,\   p=q=1$ and $2-p=2=q$,
respectively, are:
	\begin{displaymath}
P^{\pm }_{2,0}=\pm \left(
			\begin{tabular} {cc}
 0 & 1 \\
-1 & 0
			\end{tabular}
\right) , \  \
P^{\pm }_{1,1}=\pm \left(
			\begin{tabular} {cc}
 0 & 1 \\
 1 & 0
			\end{tabular}
\right) , \  \
 		P^{\pm }_{0,2}=-P^{\pm }_{2,0}.
	\end{displaymath}
For $p=n=2k$, $q=0$, we have $P=$\hbox {$diag(P_{1},\ldots ,P_{k})$},
$P_{1},\ldots ,P_{k}\in \{ P^{-}_{2,0},P^{+}_{2,0}\} $.
	\par
	In the general case, the answer to the problem
for the existence of $L$-transports consistent with a given Hermitian
structure $({\bf J},g)$ is negative.  Below we shall analyze the reasons for
this.
	\par
	Let in a fibre bundle a Hermitian structure
$({\bf J},g)$ be given. We want to see whether there exist $L$-transports along paths
consistent with it and, possibly, to describe them.
	\par First of all, for
the existence of $L$ consistent with $g$ the signature (and consequently the
number of positive eigenvalues) of $g$ must not depend on the point at which
it is (they are) calculated (see [1], proposition 2.4).
	\par
	By proposition 2.5 of [1], from the consistency between $L$ and
$g$ it follows that the matrix function $F$ describing $L$ through (1.4)
has the form (4.7) in which $Y$ is arbitrary, $Z$ satisfies (4.8)
and $p,\   q$ and $D$ are define by

	\begin{equation}
D^{\top }(x)G(x)D(x)=G_{p,q}
	\end{equation}
for any point $x$ from the base of the bundle.
\par
The consistency between ${\bf J}$ and $L$ shows that for some $C_{0}$ the
function $F$ must satisfy (4.4) and, due to (4.7), the matrix function
 $Z$ is a solution of the equation

\begin{equation}
Z^{-1}(s;\gamma )PZ(s;\gamma )=A(\gamma (s)):=
D^{-1}(\gamma (s)){\bf J}(\gamma (s))D(\gamma (s)),
\end{equation}
where $P$ is given by (4.10) for some $C_{0}$ and, in accordance with
the consistency between ${\bf J}$ and $g$ (see (4.3)), does not depend
on $\gamma $ (see the proof of proposition 4.2).
	\par
	So, the matrix function $F$ defining $L$ has the form (4.7) in
which $Z$ is a solution of the system (see (4.8) and (4.13))
		\setcounter{equation}{13}
		\def\theequation{\thechapter.\arabic{equation}a}
	\begin{equation}
		 Z^{\top }(s;\gamma )G_{p,q}Z(s;\gamma )=G_{p,q},
		\end{equation}
		\setcounter{equation}{13}
		\def\theequation{\thechapter.\arabic{equation}b}
	\begin{equation}
		 PZ(s;\gamma )-Z(s;\gamma )A(\gamma (s))=0.
	\end{equation}
		\setcounter{equation}{14}
		\def\theequation{\thechapter.\arabic{equation}}
The equalities (4.14) form a system of $n(n+1)/2+n^{2}$ scalar
equations for $n^{2}$ elements of $Z(s;\gamma )$, as a consequence of which,
generally, it has no a solution with respect to $Z$ (see, in particular,
the analysis made in [12] for the existence of solutions for
the equation $AX+XB=C$ with respect to $X$).
\par
The above consideration prove the following
	\par
	{\bf Proposition 4.3.}
	Let in a fibre bundle there be given a Hermitian
structure $({\bf J},g)$ and the signature of $g$ be independent of the point
at which it is calculated.  Let $D, \ p$ and $q$ be defined by (4.12).
Then if for every (resp. a given) path $\gamma $ and some $Y(\gamma )$
 and $C_{0}(\gamma )$ there exists a (constant) matrix $P$ satisfying (4.10)
and (4.11), for which the system (4.14) (with A defined from (4.13)) has a
solution with respect to $Z(s;\gamma )$, then the $L$-transport along paths
(resp. the given path $\gamma $) defined by the matrices (4.7) is globally
(resp. locally along $\gamma $) consistent with $({\bf J},g)$. The
$L$-transports along paths (resp. along $\gamma $) obtained
in this way form the class
of all globally (resp. locally along $\gamma $) consistent with $({\bf J},g)$
$L$-transports along paths (resp.  along $\gamma $).

\setcounter {chapter} {5}
\section {\bf CONCLUSION}
	\par In this work we have considered the problem for
consistency (or compatibility) of transports along paths in (different or
coinciding) fibre bundles and bundle morphisms between them.  Our approach
to this problem is sufficiently general and
 as its special cases includes all
known to the author analogous problems posed in the literature.  In
particular, one most often comes to the question for consistency of a
connection and some other mathematical structure, like a metric, complex or
almost complex structure.  It can equivalently be formulated as a special case
of the above problem in the following way.  On one hand, the connection
can equivalently be expressed in terms of a corresponding parallel
transport, a kind of transport along paths [13].  On the other
hand, the mentioned mathematical structures, at least, in the known to the
author analogous problems in the available to him literature, can
equivalently be put in a form of bundle morphisms of the fibre
bundle in which the parallel transport acts.  So, the consistency
between a connection and a mathematical structure is equivalent to
the consistency of a corresponding parallel transport and a bundle
morphism.  A typical example of this kind is the consistency between
a symmetric (Riemannian) metric and a linear connection (in the
tangent bundle to a manifold), which in other terms is treated by
proposition 3.2 of [1] (see also the comment after definition 4.1
of the present paper).
	\par
	In connection with proposition 2.3 there arise two problems.
First, to describe, if any, all pairs of transports (locally)
consistent along a fixed path with a given bundle morphism.  Second,
to describe, if any, all bundle morphisms (resp.  pairs of
transports along paths) globally, i.e. along every path, consistent
with a given pair of transports along paths (resp.  bundle morphism).
These problems will be investigated elsewhere.
	\par
At the end, we want to note that in the very special case when
the bundle morphism $(F,f)$ is such that there exists the inverse map
$F^{-1}$ (and hence also the map $f^{-1}$),
then all pairs of transports along paths consistent
with $(F,f)$ are $({}^{1}\! I,{}^{2}\! I)$, where
${}^{2}\! I$ is arbitrary and ${}^{1}\! I$ is given by
 ${}^{1}\! I^{\gamma }_{s \to t}=F^{-1}_{\gamma (t)}\circ
{}^{2}\! I^{f\circ \gamma }_{s \to t}\circ F_{\gamma (s)}$.
 This result is an evident corollary by eq.  (2.1).

\section* {\bf ACKNOWLEDGEMENT}
\par This research was partially supported by the Fund for Scientific
Research of Bulgaria under contract Grant No.  $F\  103.$

\section* {\bf REFERENCES}
	\begin{tabular}{rp{322pt   }}
1. & $\ \ $  Iliev B.Z., Linear transports along paths in vector bundles IV.
Consistency with bundle metrics,  E5-94-17, Dubna, 1994.\\
2. & $\ \ $  Iliev B.Z., Transports along paths in fibre bundles.  General
theory, Communication JINR, E5-93-299, Dubna, 1993.\\
3. & $\ \ $  Choquet-Bruhat Y.  at el., Analysis, manifolds and physics,
North-Holland Publ.Co., Amsterdam, 1982.\\
4. & $\ \ $  Husemoller D., Fibre bundles, McGrow-Hill Book Co., New York-St.
Louis-San Francisco-Toronto-London-Sydney, 1966.\\
5. & $\ \ $  Steenrod N., The topology of fibre bundles, 9-th ed., Princeton Univ.
Press, Princeton, 1974 (1-st ed.  1951).\\
6. & $\ \ $  Viro O.Ya., D.B. Fuks, I.  Introduction to
homotopy theory, In:  Reviews of science and technology, sec.  Modern problems
in mathematics.  Fundamental directions, vol.24, Topology-2, VINITI, Moscow,
1988, 6-121 (In Russian).\\
7. & $\ \ $  Iliev B.Z., Linear transports along paths in vector bundles.
I. General theory, Communication JINR, E5-93-239, Dubna, 1993.\\
8. & $\ \ $  Greub W., S.  Halperin, R.  Vanstone, Connections,
Curvature, and Cohomology, vol.1, vol.2, Academic Press, New York and London,
1972, 1973.\\
9. & $\ \ $  Karoubi M., K-theory.  An Introduction, Springer-Verlag,
Berlin-Heidel\-berg-New York, 1978.\\
10. & $\ \ $  Kobayashi S., K.  Nomizu, Foundations of differential geometry,
vol.  2, Interscience publishers, New-York-London, 1969.\\
11. & $\ \ $  Asanov G.  S., Finsler Geometry, Relativity and Gauge Theories,
D.  Reidel Publ.  Co., 1985; Rund H., The Deferential Geometry of
Finsler spaces, Springer Verlag, 1959.\\
12. & $\ \ $  Bellman R., Introduction to matrix analysis, McGRAW-HILL book
comp., New York-Toronto-London, 1960.\\
13. & $\ \ $  Iliev B.Z., Transports along paths in fibre bundles II.  Ties
with the theory of connections and parallel transports,
Communication JINR, E5-94-16, Dubna, 1994.\\
	\end{tabular}

\newpage
\vspace {2cm}
\Large
\indent Iliev B.  Z.
\vspace {0.41in}
\par
\begin{center}
Transports along Paths in Fibre Bundles\\
III.  Consistency with Bundle Morphisms
\end{center}
\medskip
The general problem for consistency between arbitrary
tran\-sports along paths in fibre bundles and bundle morphisms
between them is formulated and investigated.  The special case of
one fibre bundle, its morphism and transport along paths acting in it
is considered.  The consistency between linear transports
along paths in a vector bundle and a Hermitian structure in it is
studied.
\par
\medskip
\medskip
\medskip
The investigation has been performed at the Laboratory of
Theoretical Physics, JINR.
\end{document}